\documentclass[12pt,thmsa]{article}
\usepackage{amsmath, latexsym, amsfonts, amssymb, amsthm, amscd}
\usepackage{natbib}

\newtheorem{theorem}{Theorem}
\newtheorem{corollary}{Corollary}
\newtheorem{proposition}{Proposition}
\newtheorem{lemma}{Lemma}
\newtheorem{rem}{Remark}

\newcommand{\p}{\Bbb{P}}

\newcommand{\e}{\Bbb{E}}

\newcommand{\ind}{\mbox{\rm 1\hspace{-0.04in}I}}

\newcommand{\R}{\mathbb{R}}

\newcommand{{\ud}}{\textrm{d}}

\def\QED{\hfill\vrule height 1.5ex width 1.4ex depth -.1ex \vskip20pt}

\begin{document}
\hspace*{-0.5in}{\footnotesize This version oct 30, 2009.}
\vspace*{0.9in}
\begin{center}
{\Large Explicit identities for L\'evy processes associated to  symmetric stable processes.}\\
\vspace*{0.4in} {\large M.E. Caballero\footnote{$^{, 3}$ Instituto
de Matem\'aticas, Universidad Nacional Autonoma de M\'exico,
M\'exico D.F  C.P. 04510. $^1$E-mail: marie@matem.unam.mx,
$^{3}$E-mail: garmendia@matem.unam.mx}, J.C. Pardo\footnote{
Department of Mathematical Science, University of Bath. {\sc Bath}
BA2 7AY. {\sc United Kingdom},\\ E-mail:jcpm20@bath.ac.uk} and J.L.
P\'erez$^3$\vspace*{0.2in}}\\
\end{center}
\vspace{0.2in}

\begin{abstract} 
In this paper we introduce a new class of L\'evy processes which we call hypergeometric-stable L\'evy processes, because they are obtained from symmetric stable processes through several transformations and where the Gauss hypergeometric function plays an essential role. We characterize the L\'evy measure of this class and obtain several useful properties such as the Wiener Hopf factorization, the characteristic exponent and some associated exit problems.\\

\noindent {\sc Key words and phrases}: Symmetric stable L\'evy processes, Positive self-similar Markov
processes, Lamperti representation, first exit time, first hitting time.\\

\noindent MSC 2000 subject classifications: 60 G 18, 60 G 51, 60 B
52.
\end{abstract}

\vspace{0.5cm}
\section{Introduction and preliminaries.}
Let $Z=(Z_t=\{Z^{(1)}_t,\ldots Z^{(d)}_t\},t\geq 0)$ be a symmetric stable  L\'evy process of index $\alpha\in (0,2)$ in $\R^d$ ($d\ge 1$), that is, a process with stationary independent increments, its sample paths are c\`adl\`ag and
\[
\e_0\big(\exp\{i<\lambda, Z_t>\}\big)=\exp\{-t\|\lambda\|^{\alpha}\},
\]
for all $t\ge 0$ and $\lambda \in \R^d$. Here $\p_z$ denotes the law of the process $Z$ initiated  from $z\in \R^d$,  $\|\cdot\|$ the norm in  $\R^d$ and $<\cdot,\cdot>$ the Euclidean inner product.

The process $Z^{(k)}=(Z^{(k)}_t, t\ge 0)$ will be called the $k$-th coordinate process of $Z$. Of course, $Z^{(k)}$ is a real symmetric stable process whose characteristic exponent is given by
\[
\e_0\Big(\exp\big\{i\theta Z^{(k)}_t\big\}\Big)=\exp\{-t|\theta|^{\alpha}\},
\]
for all $t\ge 0$ and $\theta \in \R$.

According to Bertoin \citep{Be}, the process $Z$ is transient for $\alpha<d$, that is
\[
\lim_{t\to \infty}\|Z_t\|=\infty\qquad \textrm{a.s.},
\]
and it oscillates otherwise, i.e. for  $\alpha\in [1,2)$ and $d=1$, we have
\[
\limsup_{t\to \infty}Z_t=\infty\qquad \textrm {and} \qquad \liminf_{t\to \infty}Z_t=-\infty\qquad \textrm{a.s.}
\]
When $d\ge 2$, we have that single points are polar, i.e. for every $x,z\in \R^d$
\[
\p_x(Z_t=z\quad \textrm{for some }t>0)=0.
\]
In the one-dimensional case, points are polar for $\alpha\in(0,1]$ and when $\alpha\in(1,2)$  the process $Z$ makes infinitely many jumps across a point, say $z$, before the first hitting time of $z$ (see for instance Proposition VIII.8 in \citep{Be}).

One of the main properties of the process  $Z$ is that it  satisfies the scaling property with index $\alpha$, i.e. for every $b>0$
\begin{equation}\label{scale}
\mbox{\it \textrm{The law of} $\;(bZ_{b^{-\alpha}t},\,t\ge0)$ \textrm{under} $\p_x$ \textrm{is}
$\p_{bx}$.}
\end{equation}
This implies that  the radial process $R=(R_t, t\ge 0)$ defined by  $R_t=\|Z_t\|$ satisfies the same scaling property (\ref{scale}). Since $Z$ is isotropic, its radial part  $R$ is a strong Markov process (see Millar \citep{mi}).  When $d\ge 2$, the radial process $R$ hits points if and only if $Z^{(1)}$  hits points i.e. when $\alpha\in(1,2)$ (see for instance Theorem 3.1 in \citep{mi}). Finally, we note that when  points are polar for $Z$  the radial process $R$ will never hit the point $0$.

In what follows we will assume that $\alpha\le d$, so the  radial process $R$ will be  a positive self-similar Markov process (pssMp) with index $\alpha$ and infinite lifetime.
 A natural question arises: 
 can we characterize the L\'evy process $\xi$ associated to the pssMp  $(R_t, t\ge 0)$ via the 
 Lamperti transformation?

We briefly recall  the main features of the Lamperti transfomation, between pssMp and L\'evy processes.    A positive  self-similar
Markov processes $(X,\mathbb{Q}_x)$, $x>0$, is a strong Markov processes with
c\`adl\`ag paths, which fulfills a scaling property.   Well-known examples of this kind of processes are: Bessel processes, stable subordinators, stable processes conditioned to stay positive, etc.

 According to Lamperti
\citep
{La}, any pssMp up to its first hitting time of 0 may be
expressed as the exponential of a L\'evy process, time changed by
the inverse of its exponential functional. More formally, let
$(X,\mathbb{Q}_x)$ be a pssMp with index $\beta>0$, starting from $x>0$,
set
\[S = \inf \{t>0 : X_t =0\}\]
and write the canonical process $X$ in the following form:
\begin{equation}\label{lamp}
X_t=x\exp\left\{\xi_{\tau(tx^{-\beta})}\right\} \qquad 0\le t<S\,,
\end{equation}
where for $t<S$,
\[\tau(t) = \inf \left\{s\geq 0 : \int_0^s
\exp\left\{\beta\xi_u\right\} {\rm d} u \geq t\right\}.\] Then under
$\mathbb{Q}_x$, $\xi=(\xi_t,\;t\geq 0)$ is a L\'evy process started from $0$
whose law does not depend on $x>0$ and such that:
\begin{itemize}
\item[$(i)$] if  $\mathbb{Q}_x(S=+\infty)=1$,  then  $\xi$ has an infinite lifetime
and $\displaystyle\limsup_{t\rightarrow+\infty}\xi_t=+\infty$,
$\p_x$-a.s.,

\item[$(ii)$] if $\mathbb{Q}_x(S<+\infty,\,X(S-)=0)=1$, then $\xi$ has an infinite
lifetime and $\displaystyle\lim_{t\to\infty} \xi_t = -\infty$,
$\p_x$-a.s.,

\item[$(iii)$] if $\mathbb{Q}_x(S<+\infty,\,X(S-)>0)=1$,  then $\xi$ is killed at an
independent exponentially distributed random time with parameter
$\lambda>0$.

\end{itemize}
As mentioned in \citep
{La}, the probabilities
$\mathbb{Q}_x(S=+\infty)$, $\mathbb{Q}_x(S<+\infty,\,X(S-)=0)$ and
$\mathbb{Q}_x(S<+\infty,\,X(S-)>0)$ are 0 or 1 independently of $x$, so that
the three classes presented above are exhaustive. Moreover, for any
$t < \int_0^{\infty}\exp\{\beta\xi_s\}\,{\rm d} s$,
\begin{equation}\label{1664}\tau(t)=\int_0^{x^\beta t}\frac{{\rm d} s}
{(X_s)^\beta}\,,\;\;\; \mathbb{Q}_x-\mbox{a.s.}\end{equation} Therefore
(\ref{lamp}) is invertible and yields a one-to-one relation between
the class of pssMp's killed at time $S$ and the one of L\'evy
processes.

Another important result of  Lamperti \citep
{La} provides the explicit form of the generator of any pssMp $(X, \mathbb{Q}_y)$ in terms of its underlying L\'evy process.  Let $\xi$ be the underlying L\'evy process associated to $(X, \mathbb{Q}_y)$ via (\ref{lamp}) and denote by $\mathcal{L}$ and $\mathcal{M}$  their respective infinitesimal generators.  Let $\mathcal{D}_{\mathcal{L}}$ be the domain of the generator $\mathcal{L}$ and recall that it contains all the functions with continuous second derivatives on $[-\infty, \infty]$, and that if $f$ is such a function then $\mathcal{L}$ acts as follows for $x\in \R$, where $\mu\in\R$ and $\sigma>0$:
\begin{equation}\label{igLp}
\mathcal{L}f(x)=\mu f'(x)+\frac{\sigma^2}{2}f''(x)+\int_{\footnotesize\R}\big(f(x+y)-f(x)-f'(x)\ell(y)\big)\Pi({\rm d} y)-bf(x).
\end{equation}
The measure $\Pi({\rm d} x)$ is the so-called L\'evy measure of $\xi$, which  satisfies
\[
\Pi(\{0\})=0\qquad \textrm{ and } \qquad\int_{\footnotesize\R} (1\land |x|^2)\Pi({\rm d} x)<\infty.
\]
The function $\ell(\cdot)$ is a bounded Borel function such that $\ell(y)\sim y$ as $y\to 0$. The positive constant $b$  represents the killing rate of $\xi$ (b=0 if $\xi$ has infinite lifetime).
Lamperti  establishes the following result in \citep
{La}.
\begin{theorem}\label{thmlamp}
If $g$ is such that $g$, $yg'$ and $y^2g''$ are continuous on $[0,\infty]$, then they belong to the domain, $\mathcal{D}_\mathcal{M}$, of the infinitesimal generator of  $(X, \mathbb{Q}_y)$, which acts as follows for $y>0$
\[
\begin{split}
\mathcal{M}g(y)=\mu y^{1-\beta}g'(y)&+\frac{\sigma^2}{2}y^{2-\beta}g''(y)-by^{-\beta}g(y)\\
&+y^{-\beta}\int_{0}^{\infty}\big(g(yu)-g(y)-g'(y)\ell(\log u)\big)G({\rm d}u),
\end{split}
\]
where $G({\rm d} u)=\Pi({\rm d} u)\circ\log u $, for $u>0$. This expression determines the law of the process $(X_t, 0\le t\le T)$ under $\mathbb{Q}_y$.
\end{theorem}

 Previous work on this subject  appears in  Carmona et al. \citep
{CPY} where the authors
 studying the radial part of  a Cauchy process $C=(C_t, t\geq 0)$ (i.e. $\alpha=d=1$), they obtain the infinitesimal generator of its associated  L\'evy process $\xi=(\xi_t, t\geq 0)$ via the Lamperti transformation. More preciseley, the infintesimal generator of   $\xi$ is given as follows
 $$\mathcal{L}g(\xi)=\frac{1}{\pi}\int\frac{\cosh \eta}{(\sinh \eta)^2}(g(\xi+\eta)-g(\xi)-\eta g^{\prime}(\xi \ind_{\vert\eta\vert \leq 1})d\eta,$$
and  its  characteristic exponent satisfies
  $$\e\Big(\exp\{i\lambda \xi_t\}\Big)= e^{-i\lambda \tanh\frac{\pi\lambda}{2}}.$$
As we will see in sections 2 and 5 this example is a particular case of the results obtained in this paper by very different methods. As it is expected, the formulas obtained in both papers coincide  for $\alpha=d=1$.
  
It is important to point out that in Carmona et al. \citep
{CPY}, it is announced that the authors will continue this line of reseach by studying the case of the norm of a multidimensional  Cauchy process, but up to our knowledge this has not be done.
   
 The paper is organized as follows: In section 2, we compute the infinitesimal generator of the radial process R and using theorem
  \ref{thmlamp} we obtain the characteristics of its associated L\'evy process $\xi$. The L\'evy measure obtained has a rather complicated form since it is expressed in terms of the Gauss hypergeometric function $_{2}\mathcal{F}_{1}$. When $d=1$ we show that the process $\xi$ can be expressed as the sum of a Lamperti stable process (see Caballero et al.\citep{cpp} for a proper definition) and an independent  Poisson process. 
    
  In section 3 we study  one sided exit problems of the L\'evy process $\xi$, using well known results of Blumenthal et al. \citep
{Bl} for the symmetric $\alpha$-stable process $Z$. When $\alpha<d$, a straightforward computations allows us to deduce the law of the random variable $\underline{\xi}_{\infty}=\inf_{t\ge 0}\xi_t$.
  
  In section 4,  we study the special  case  $1<\alpha<d$. Using the work of S. Port \citep
{po} on the radial processes of $Z$, we compute the probability that the L\'evy process $\xi$ hits points.
  
    Finally in section 5  we obtain the Wiener-Hopf factorization of $\xi$ and deduce the explicit form of the characteristic exponent. 
    Concluding remarks show in section 6  how to obtain n-tuple laws  for $\xi$ and $R$ following Kyprianou et al. \citep
{kpr}.

\section{The underlying L\'evy process of $R$}
In this section, we compute the generator of the radial process $R$ and the  characteristics of the underlying L\'evy process $\xi$ in the Lamperti representation (\ref{lamp}) of the latter.

To this end, it will be useful to invoke the expression of $Z$ as a subordinated Brownian motion. More precisely, let $B=(B_t,t\ge 0)$ be a $d$-dimensional Brownian motion initiated from $x\in \R^d$ and let $\sigma=(\sigma_t,t\ge 0)$ be an independent stable subordinator with index $\alpha/2$ initiated from $0$. Then the process $(B_{2\sigma_t}, t\ge 0)$ is a standard symmetric  $\alpha$-stable process.

Let us define the so-called  Pochhammer symbol by
\[
(z)_\alpha=\frac{\Gamma(z+\alpha)}{\Gamma(z)}, \qquad \textrm{ for }\quad z\in \mathbb{C},
\]
 and the Gauss's hypergeometric function by
\[
_{2}\mathcal{F}_{1}\Big(a,b;c;z\Big)=\sum_{k=0}^{\infty}z^{k}\frac{(a)_{k}(b)_{k}}{(c)_{k}\,k!}, \qquad \textrm{for }\quad \|z\|<1,
\]
where $a,b,c>0$.
\begin{theorem}
If $g:\R_+\to \R$ is such that $g\in C^2_0(\R_+)$. Hence the infinitesimal generator of $R=(R_t,t\ge 0)$, denoted by $M$, acts  as follows for $a>0$,
\[
\begin{split}
Mg(a)&= a^{-\alpha}\int_{0}^{\infty}\big(g(y a)-g(a)-g'(a)\ell(\log y)\big)\frac{y^{d-1}}{(1+y^{2})^{(\alpha+d)/2}}\overline{F}\left(\left(\frac{2y}{1+y^{2}}\right)^{2}\right){\rm d} y,\notag\\
\end{split}
\]
where
\begin{equation}\label{gausshyp}
\overline{F}(z)=\frac{2^{\alpha}\alpha(d/2)_{\alpha/2}}{\Gamma(1-\alpha/2)}\, _{2}\mathcal{F}_{1}\Big((\alpha+d)/4,(\alpha+d)/4+1/2;d/2;z\Big) \qquad \textrm{for }\quad z\in(-1,1),
\end{equation}
and the
function $\ell$ is given by
\begin{equation}\label{funl}
\ell(y)=\frac{y}{1+y^{2}}e^{(1-d)y}\big(1+e^{2y}\big)^{(\alpha+d)/2-1}\ind_{\{|y|<1\}}.
\end{equation}
\end{theorem}
\noindent{\it Proof:} From Theorem 32.1 in \citep
{Sa} and the fact that $Z$ can be seen as a subordinated Brownian motion, the
infinitesimal generator $M$ of $R=(R_t,t\ge 0)$ is given as follows
\begin{equation}\label{g:1}
Mh=\int_{0}^{\infty}(P_{s}h-h)\rho({\rm d} s),
\end{equation}
where  $\rho$ is the L\'evy measure of the stable subordinator $2\sigma$ and is given by
\[
\rho({\rm d} s)=\frac{2^{\alpha/2-1}\alpha}{\Gamma(1-\alpha/2)}s^{-(1+\alpha/2)}\ind_{\{s>0\}}{\rm d} s
,\]
 $P_{s}$ is the  semi-group of the $d$-dimensional Bessel process and $h$ is any function in the domain of the infinitesimal generator of $(P_t,t\ge 0)$.

Let $g$ be  as in the statement and  recall that for $a>0$, the semi-group for the
$d$-dimensional Bessel process satisfies
\[
P_{s}g(a)=\int_{0}^{\infty}{\rm d} y\, \frac{g(y)}{s}\left(\frac{y}{a}\right)^{d/2-1}y
\exp\left(-\frac{y^2+a^2}{2s}\right)I_{d/2-1}\left(\frac{a y}{s}\right),
\]
where $I_{d/2-1}$ is the modified Bessel function of index $d/2-1$
(see for instance \citep
{ry}). Therefore putting the pieces together,
it follows
\begin{align}\label{gen}
Mg(a)&=\frac{2^{\alpha/2-1}\alpha}{\Gamma(1-\alpha/2)}\int_{0}^{\infty}\int_{0}^{\infty}y\Big(g(y)-g(a)\Big)\left(\frac{y}{a}\right)^{d/2-1}\notag\\
&\hspace{6cm}\times\frac{1}{s^{2+\alpha/2}}
\exp\left(-\frac{a^2+y^{2}}{2s}\right)I_{d/2-1}\left(\frac{ay}{s}\right){\rm
d} y{\rm d} s.
\end{align}
Now, recall the following identity of the modified Bessel function
$I_{d/2-1}$,
\begin{equation}
I_{d/2-1}(x)=\sum_{k=0}^{\infty}\frac{(x/2)^{2k+d/2-1}}{\Gamma(d/2+k)k!},\notag
\end{equation}
and note that for $a\ne y$
\begin{align}\label{a:4}
&\hspace{-.8cm}\int_{0}^{\infty}\frac{{\rm d} s}{s^{2+\alpha/2}}
\exp\left(-\frac{a^2+y^{2}}{2s}\right)I_{d/2-1}\left(\frac{ay}{s}\right)\notag\\
&=\sum_{k=0}^{\infty}\int_{0}^{\infty}{\rm d} s\left(\frac{ay}{2s}\right)^{2k+d/2-1}\frac{s^{-2-\alpha/2}}{\Gamma(d/2+k)k!}\exp\left(-\frac{a^{2}+y^{2}}{2s}\right)\notag\\
&=\sum_{k=0}^{\infty}\frac{1}{\Gamma(d/2+k)k!}\left(\frac{ay}{\alpha^2+y^2}\right)^{2k+(\alpha+d)/2}\left(\frac{2}{ay}\right)^{1+\alpha/2}\int_{0}^{\infty}{\rm d} u\,u^{2k+(\alpha+d)/2-1}e^{-u}\notag\\
&=2^{1+\alpha/2}\frac{(ay)^{d/2-1}}{(a^{2}+y^{2})^{(\alpha+d)/2}}\sum_{k=0}^{\infty}\left(\frac{ay}{a^{2}+y^{2}}\right)^{2k}\frac{\Gamma(2k+(\alpha+d)/2)}{\Gamma(k+1)\Gamma(d/2+k)}.
\end{align}
Next, we consider the following property of the Gamma function,
\begin{equation}\label{propgf}
\Gamma(2z)=(2\pi)^{-1/2}2^{2z-1/2}\Gamma(z)\Gamma(z+1/2),
\end{equation}
and deduce that
\begin{align}
\Gamma(2k+(\alpha+d)/2)&=
(2\pi)^{-1/2}2^{2k+(\alpha+d)/2-1/2}\Gamma(k+(\alpha+d)/4)\Gamma(k+(\alpha+d)/4+1/2)\notag\\
&=2^{2k}\Gamma((\alpha+d)/2)((\alpha+d)/4)_{k}((\alpha+d)/4+1/2)_{k}.\notag
\end{align}
Therefore using the above identity,  we see that (\ref{a:4}) is equal to
\begin{align}
\frac{2^{\alpha/2+1}(ay)^{d/2-1}}{(a^{2}+y^{2})^{(\alpha+d)/2}}\frac{\Gamma((\alpha+d)/2)}{\Gamma(d/2)}\sum_{k=0}^{\infty}\left(\left(\frac{2ay}{a^{2}+y^{2}}\right)^{2}\right)^{k}\frac{((\alpha+d)/4)_{k}((\alpha+d)/4+1/2)_{k}}{(d/2)_{k}\,k!},\notag
\end{align}
where the series above is the Gauss's hypergeometric function
\begin{align*}
_{2}\mathcal{F}_{1}&\left((\alpha+d)/4,(\alpha+d)/4+1/2;d/2;\left(\frac{2ay}{a^{2}+y^{2}}\right)^{2}\right).
\end{align*}
We remark that we cannot use Fubini's theorem on (\ref{gen}) because
the expression inside the integral with respect to the product
measure is not integrable. This is easily seen by noting that
\begin{equation}
\left|_{2}\mathcal{F}_{1}\left((\alpha+d)/4,(\alpha+d)/4+1/2;d/2;\left(\frac{2ay}{a^{2}+y^{2}}\right)^{2}\right)\right|\sim|y-a|^{-(\alpha+1)}\qquad\text{as
$y\to a$}.\notag
\end{equation}
So instead let us consider $\varepsilon_1,\varepsilon_2$, $c\geq0$,
and denote by
$$A_{\varepsilon_1,\varepsilon_2}(c)=\{y\in(0,\infty):y>c+\varepsilon_1\}\cup\{y\in(0,\infty):y<c-\varepsilon_2/(c+\varepsilon_1)\}.$$
Then we have
\begin{align}\label{arreglo1}
&\int_{0}^{\infty}\int_{A_{\varepsilon,a\varepsilon}(a)}y\Big(g(y)-g(a)\Big)\left(\frac{y}{a}\right)^{d/2-1}\frac{1}{s^{2+\alpha/2}}
\exp\left(-\frac{a^2+y^{2}}{2s}\right)I_{d/2-1}\left(\frac{ay}{s}\right){\rm
d} y{\rm d} s.
\end{align}
We would like to use Fubini's Theorem in the expression above, to
this end we now prove the integrability of the integrand with
respect the product measure. For simplicity, we use the notation
established in (\ref{gausshyp}), and using Tonelli's theorem and
(\ref{a:4}) we have
\begin{align*}
&\int_{0}^{\infty}\int_{A_{\varepsilon,a\varepsilon}(a)}y\Big|g(y)-g(a)\Big|\left(\frac{y}{a}\right)^{d/2-1}\frac{1}{s^{2+\alpha/2}}
\exp\left(-\frac{a^2+y^{2}}{2s}\right)I_{d/2-1}\left(\frac{ay}{s}\right){\rm
d} y{\rm d} s.\\
&\leq
2\|g\|_{\infty}\int_{A_{\varepsilon,a\varepsilon}(a)}\frac{y^{d-1}}{(a^{2}+y^{2})^{(\alpha+d)/2}}\overline{F}\left(\left(\frac{2ay}{a^{2}+y^{2}}\right)^{2}\right){\rm
d} y,
\end{align*}
which is finite.
 So now let us return to (\ref{arreglo1}), then
applying Fubini's theorem and (\ref{a:4}) we obtain
\begin{align}\label{arreglo2}
&\int_{0}^{\infty}\int_{A_{\varepsilon,a\varepsilon}(a)}y\Big(g(y)-g(a)\Big)\left(\frac{y}{a}\right)^{d/2-1}\frac{1}{s^{2+\alpha/2}}
\exp\left(-\frac{a^2+y^{2}}{2s}\right)I_{d/2-1}\left(\frac{ay}{s}\right){\rm
d} y{\rm d} s.\notag\\
&=\int_{A_{\varepsilon,a\varepsilon}(a)}\Big(g(y)-g(a)\Big)\frac{y^{d-1}}{(a^{2}+y^{2})^{(\alpha+d)/2}}\overline{F}\left(\left(\frac{2ay}{a^{2}+y^{2}}\right)^{2}\right){\rm
d} y\notag\\
&=a^{-\alpha}\int_{C(a,\varepsilon)}\Big(g(ay)-g(a)\Big)\frac{y^{d-1}}{(1+y^{2})^{(\alpha+d)/2}}\overline{F}\left(\left(\frac{2y}{1+y^{2}}\right)^{2}\right){\rm
d} y.
\end{align}
 where $C(a,\varepsilon)=\{y: 0<y<\frac{a}{a+\varepsilon}\}\cup\{y: 1+\frac{\varepsilon}{a}<y\}.$
In order to get the result, we first show that  if
$$B(a,\varepsilon)=\left(\frac{1}{e},\frac{a}{a+\varepsilon}\right)\cup\left(1+\frac{\varepsilon}{a},e\right)=C(a,\varepsilon)\cap(1/e,e),$$ then
\begin{equation}\label{inte0}
\int_{B(a,\varepsilon)} \frac{\log
y}{1+\log^{2}y}\frac{1}{1+y^{2}}\overline{F}\left(\left(\frac{2y}{1+y^{2}}\right)^{2}\right){\rm
d} y=0.
\end{equation}
To do so, we note that the integral in (\ref{inte0}) is equal to
\[
\int_{1+a^{-1}\varepsilon}^{e}\frac{\log
y}{1+\log^{2}y}\frac{1}{1+y^{2}}\overline{F}\left(\left(\frac{2y}{1+y^{2}}\right)^{2}\right){\rm
d} y +\int_{1/e}^{a/(a+\varepsilon)}\frac{\log
y}{1+\log^{2}y}\frac{1}{1+y^{2}}\overline{F}\left(\left(\frac{2y}{1+y^{2}}\right)^{2}\right){\rm
d} y.
\]
Making the change of variable $y=z^{-1}$ in the first integral of above, we get that
\begin{align}
\int_{1+a^{-1}\varepsilon}^{e}\frac{\log y}{1+\log^{2}y}\frac{1}{1+y^{2}}&\overline{F}\left(\left(\frac{2y}{1+y^{2}}\right)^{2}\right){\rm d} y\notag\\
&=-\int_{1/e}^{a/(a+\varepsilon)}\frac{\log
z}{1+\log^{2}z}\frac{1}{1+z^{2}}\overline{F}\left(\left(\frac{2z}{1+z^{2}}\right)^{2}\right){\rm
d} z,\notag
\end{align}
and the identity (\ref{inte0}) follows. It is easy to see using
(\ref{a:4}) the following equality:
\begin{align}\label{arreglo3a}
&\int_{B(a,\varepsilon)}\frac{\log
y}{1+\log^{2}y}\frac{1}{1+y^{2}}\overline{F}\left(\left(\frac{2y}{1+y^{2}}\right)^{2}\right){\rm
d} y\notag\\
&=\frac{a^{\alpha}2^{\alpha/2-1}\alpha}{\Gamma(1-\alpha/2)}\int_{0}^{\infty}\int_{0}^{\infty}y\ell(\log y/a)\left(\frac{y}{a}\right)^{d/2-1}\ind_{B(a,\varepsilon)}(y)\notag\\
&\hspace{6cm}\times\frac{1}{s^{2+\alpha/2}}
\exp\left(-\frac{a^2+y^{2}}{2s}\right)I_{d/2-1}\left(\frac{ay}{s}\right){\rm
d} y{\rm d} s.
\end{align}
where $\ell$ is defined as in (\ref{funl}). Finally, we add the term
\[
a^{-\alpha}\int_{0}^{\infty}g'(a)\frac{\log
y}{1+\log^{2}y}\frac{1}{1+y^{2}}\overline{F}\left(\left(\frac{2y}{1+y^{2}}\right)^{2}\right)\ind_{B(a,\varepsilon)}(y){\rm
d} y,
\]
to the identity (\ref{arreglo2}) and after some
calculations using (\ref{arreglo3a}) we obtain
\begin{align}\label{arreglo3}
&\frac{2^{\alpha/2-1}\alpha}{\Gamma(1-\alpha/2)}\int_{0}^{\infty}\int_{A_{\varepsilon,a\varepsilon}(a)}y\Big(g(y)-g(a)-g'(a)\ell(\log
(y/a))\Big)\left(\frac{y}{a}\right)^{d/2-1}\notag\\
&\hspace{6cm}\times\frac{1}{s^{2+\alpha/2}}
\exp\left(-\frac{a^2+y^{2}}{2s}\right)I_{d/2-1}\left(\frac{ay}{s}\right){\rm
d} y{\rm d} s\notag\\
&=a^{-\alpha}\int_{B(a,\varepsilon)}\big(g(ya)-g(a)-g'(a)\ell(\log
y)\big)\frac{y^{d-1}}{(1+y^{2})^{(\alpha+d)/2}}\overline{F}\left(\left(\frac{2y}{1+y^{2}}\right)^{2}\right){\rm
d} y.
\end{align}
So using the dominated convergence theorem and (\ref{arreglo3}), we
can conclude that
\begin{align*}
Mg(a)&=\frac{2^{\alpha/2-1}\alpha}{\Gamma(1-\alpha/2)}\int_{0}^{\infty}\int_{0}^{\infty}y\Big(g(y)-g(a)\Big)\left(\frac{y}{a}\right)^{d/2-1}\notag\\
&\hspace{6cm}\times\frac{1}{s^{2+\alpha/2}}
\exp\left(-\frac{a^2+y^{2}}{2s}\right)I_{d/2-1}\left(\frac{ay}{s}\right){\rm
d} y{\rm d} s\notag\\
&=\lim_{\varepsilon \to
0}\frac{2^{\alpha/2-1}\alpha}{\Gamma(1-\alpha/2)}\int_{0}^{\infty}\int_{A_{\varepsilon,a\varepsilon}(a)}y\Big(g(y)-g(a)-g'(a)\ell(\log
(y/a))\Big)\left(\frac{y}{a}\right)^{d/2-1}\notag\\
&\hspace{6cm}\times\frac{1}{s^{2+\alpha/2}}
\exp\left(-\frac{a^2+y^{2}}{2s}\right)I_{d/2-1}\left(\frac{ay}{s}\right){\rm
d} y{\rm d} s\notag\\
&=a^{-\alpha}\int_{0}^{\infty}\big(g(ya)-g(a)-g'(a)\ell(\log
y)\big)\frac{y^{d-1}}{(1+y^{2})^{(\alpha+d)/2}}\overline{F}\left(\left(\frac{2y}{1+y^{2}}\right)^{2}\right){\rm
d} y.
\end{align*}
\QED
Using Lamperti's result (recalled in Theorem \ref{thmlamp}) and
Proposition 1, we may now give the explicit form of the generator of $\xi$. We will call this new class of L\'evy processes {\bf hypergeometric-stable}. 
\begin{corollary}
Let $\xi$ be the L\'evy process in the Lamperti representation (\ref{lamp}) of the radial process $R$. The infinitesimal generator $\mathcal{A}$, of $\xi$, with domain $\mathcal{D}_{\mathcal{A}}$ is given in the
polar case
\begin{equation}\label{LpRp}
\mathcal{A}f(x)=\int_{\footnotesize\R}\big(f(x+y)-f(x)-f'(x)\ell(y)\big)\Pi({\rm d} y),
\end{equation}
for any $f\in \mathcal{D}_{\mathcal{A}}$ and $x\in \R$, where
\[
\Pi({\rm d}y)=\frac{e^{dy}}{(1+e^{2y})^{(\alpha+d)/2}}\overline{F}\left(\frac{4e^{2y}}{(e^{2y}+1)^{2}}\right){\rm d} y.
\]
Equivalently,  the characteristic exponent of $\xi$ is given by
\[
\Psi(\lambda)=i\lambda b+\int_{\footnotesize\R}\Big(1-e^{i\lambda y}+i\lambda y\ind_{\{|y|<1\}}\Big)\Pi({\rm d} y)
\]
where
\begin{align}
b=\int_{\mathbb{R}}\Big(\ell(y)-y\ind_{\{|y|\leq1\}}\Big)\frac{e^{d y}}{(1+e^{2y})^{(\alpha+d)/2}}\overline{F}\left(\frac{4e^{2y}}{(e^{2y}+1)^{2}}\right){\rm d} y\notag.
\end{align}
\end{corollary}
We finish this section with a remarkable result on the decomposition of the L\'evy measure of the  process $\xi$ when the dimension is $d=1$ and $\alpha \in (0,1]$ (polar case). Such decomposition describes the structure of $\xi$ in terms of two independent  L\'evy processes, each with different
type of path behaviour.

  Recall in this case that the symmetric stable process $Z$ is of bounded variation  and so its radial part  $R$ and the L\'evy process $\xi$.  Hence, the characteristic exponent of $\xi$ is given by
\[
\Psi(\lambda)=\int_{\footnotesize\R}\Big(e^{i\lambda y}-1\Big)\Pi({\rm d} y).
\]
\begin{proposition}
Assume that $d=1$, then we have
\[
\Psi(\lambda)=\int_{\footnotesize\R}\Big(e^{i\lambda y}-1\Big)\Pi_1({\rm d} y)+\int_{\footnotesize\R}\Big(e^{i\lambda y}-1\Big)\Pi_2({\rm d} y),
\]
where $\Pi_{1}$ is the L\'evy measure of a Lamperti L\'evy process with characteristics $(0,1,\alpha)$ (see for instance \citep
{cpp}),  i.e.
\[
\Pi_{1}({\rm d} y)=\frac{2^{\alpha-1}\alpha(1/2)_{\alpha/2}}{\Gamma(1-\alpha/2)}\left(\frac{e^{y}}{(e^{y}-1)^{\alpha+1}}1_{\{y>0\}}+\frac{e^{y}}{(1-e^{y})^{\alpha+1}}1_{\{y<0\}}\right){\rm d} y,
\]
and
\[
\Pi_{2}({\rm d} y)=\frac{2^{\alpha-1}\alpha(1/2)_{\alpha/2}}{\Gamma(1-\alpha/2)}\frac{e^{y}}{(e^{y}+1)^{\alpha+1}}{\rm d} y,
\]
is the L\'evy measure of a compound Poisson process.
\end{proposition}
\noindent{\it Proof:} Let $x\in[0,1)$. Using identity (\ref{propgf}) twice, we deduce
\begin{align*}
_{2}\mathcal{F}&_{1}\Big((\alpha+1)/4,(\alpha+1)/4+1/2;1/2;x^{2}\Big)=\sum_{k=0}^{\infty}x^{2k}\frac{((\alpha+1)/4)_{k}((\alpha+1)/4+1/2)_{k}}{k!(1/2)_{k}}\notag\\
&\hspace{2cm}=\frac{\Gamma(1/2)}{\Gamma((\alpha+1)/4+1/2)}\frac{2^{1/2-\alpha/2}}{\Gamma((\alpha+1)/4)}\sum_{k=0}^{\infty}x^{2k}\frac{\Gamma((\alpha+1)/2+2k)}{\Gamma(2k+1)}\notag\\
&\hspace{2cm}=\frac{2^{1/2-\alpha/2}\Gamma(1/2)}{(2\pi)^{1/2}2^{1/2-(\alpha+1)/2}\Gamma((\alpha+1)/2)}\\
&\hspace{4cm}\times\frac{1}{2}\left(\sum_{0}^{\infty}x^{k}\frac{\Gamma((\alpha+1)/2+k)}{\Gamma(1+k)}+\sum_{0}^{\infty}(-x)^{k}\frac{\Gamma((\alpha+1)/2+k)}{\Gamma(1+k)}\right)\notag\\
&\hspace{2cm}=\frac{1}{2}\left(\sum_{k=0}^{\infty}x^{k}\frac{((\alpha+1)/2)_{k}}{k!}+\sum_{k=0}^{\infty}(-x)^{k}\frac{((\alpha+1)/2)_{k}}{k!}\right)\notag\\
&\hspace{2cm}=2^{-1}\Big((1-x)^{-(\alpha+1)/2}+(1+x)^{-(\alpha+1)/2}\Big).
\end{align*}
Now, from the above identity, we deduce that the L\'evy
measure of the process $\xi$ satisfies
\begin{align}
\Pi({\rm d} y)&=\frac{2^{\alpha-1}\alpha(1/2)_{\alpha/2}}{\Gamma(1-\alpha/2)}\frac{e^{y}}{(1+e^{2y})^{(\alpha+1)/2}}\left(\left(1-\frac{2e^{y}}{e^{2y}+1}\right)^{-\frac{\alpha+1}{2}}+\left(1+\frac{2e^{y}}{e^{2y}+1}\right)^{-\frac{\alpha+1}{2}}\right){\rm d} y\notag\\
&=\frac{2^{\alpha-1}\alpha(1/2)_{\alpha/2}}{\Gamma(1-\alpha/2)}e^{y}\left(\frac{1}{|e^{y}-1|^{\alpha+1}}+\frac{1}{(e^{y}+1)^{\alpha+1}}\right){\rm d} y,\notag
\end{align}
and the statement follows.\QED
\section{Entrance laws for the process $\xi$: Intervals.}
In this section, we study the probability that  the hypergeometric-stable L\'evy
process $\xi$ makes its first exit from an interval. In particular,  we obtain some explicit identities for the  one-sided exit problems.

In what follows, $P$ will be a reference probability measure on $\mathcal{D}$ (the Skorokhod space of $\R$-valued c\`adl\`ag paths) under which $\xi$ is the hypergeometric-stable L\'evy process described in Corollary 1 starting from $0$.  For any $y\in\mathbb{R}$ let
\[
T^+_y=\inf\{t\ge 0:\xi_t> y\}\;\;\mbox{and}\;\;T_y^-=\inf\{t\ge 0:\xi_t<
y\}\,,
\]
and for any $x>0$ let
\[
 \sigma^+_x=\inf\{t\ge 0:R_t>
x\}\;\;\mbox{and}\;\;\sigma_x^-=\inf\{t\ge 0:R_t < x\}.
\]

\begin{lemma}\label{generic}Fix $-\infty< v<0<u<\infty$.
Suppose that $A$ is any interval in $[u,\infty)$ and $B$ is any
interval in $(-\infty, v]$. Then,
\[
 P\Big(\xi_{T^+_u}\in A; T^+_u< \infty\Big) =
 \mathbb{P}_x\Big(R_{\sigma^+_{e^u}} \in e^A; \sigma^+_{e^u}< \infty\Big)
\]
and
\[
P\Big(\xi_{T^-_v}\in B;  T^-_v<\infty\Big) =
\mathbb{P}_x\Big(R_{\sigma^-_{e^v}} \in e^B;
\sigma^-_{e^v}<\infty\Big),
\]
where $x$ satisfies that $\|x\|=1$.
\end{lemma}
The proof is a consequence of the Lamperti
representation  and is left as an exercise. Although somewhat obvious, this lemma indicates that in order to understand the exit problem for the process $\xi$, we need to study how the radial process $R$ exits a positive interval around $x>0$. Fortunately this is possible thanks to a result of Blumenthal et al. \citep
{Bl} who established the following for the symmetric  $\alpha$-stable process $Z$.

Define,
\[
f(y,z)=\pi^{-(d/2+1)}\Gamma\left(\frac{d}{2}\right)\sin\left(\frac{\pi\alpha}{2}\right)\big|1-\|y\|^{2}\big|^{\alpha/2}\big|1-\|z\|^{2}\big|^{-\alpha/2}\|y-z\|^{-d}.
\]
\begin{theorem}[Blumenthal et al. \citep
{Bl}]Suppose that $\alpha<d$ and that $(Z,\p_x)$ is a symmetric $\alpha$-stable process with values in $\R^d$, initiated from $x$. For $\|y\|<1$ and $\|z\|\ge1$, we have
\begin{align}\label{bl}
\p_{y}\Big(Z_{\sigma^{+}_{1}}\in {\rm d} z; \sigma^{+}_{1}<\infty\Big)=f(y,z){\rm d} z.
\end{align}
Similarly for $\|y\|>1$ and $\|z\|\le1$, we have
\begin{align}\label{bla}
\p_{y}\Big(Z_{\sigma^{-}_{1}}\in {\rm d} z; \sigma^{-}_{1}<\infty\Big)=f(y,z){\rm d} z.
\end{align}
\end{theorem}
The one-side exit problem for $\xi$ can be solved using  Lemma 1 and Theorem 3 as follows.
\begin{theorem} Suppose that $\alpha<d$ and
fix $\theta\geq 0$ and $-\infty<v<0<u<\infty$. Then
\begin{align}\label{teo4.1}
P&\left(\xi_{T_{u}^{+}}-u\in
{\rm d}\theta,T_{u}^{+}<\infty\right)\notag\\
&\hspace{1.5cm}=\frac{2}{\pi}\sin\left(\frac{\pi\alpha}{2}\right)e^{2(u+\theta)}\big(1-e^{-2u}\big)^{\alpha/2}\big(e^{2\theta}-1\big)^{-\alpha/2}\big(e^{2(\theta+u)}-1\big)^{-1}{\rm d}\theta,
\end{align} and
\begin{align}\label{teo4.2}
P&\left(v-\xi_{T_{v}^{-}}\in
{\rm d}\theta,T_{v}^{-}<\infty\right)\notag\\
&\hspace{1.5cm}=\frac{2}{\pi}\sin\left(\frac{\pi\alpha}{2}\right)e^{d(v-\theta)}\big(e^{-2v}-1\big)^{\alpha/2}\big(1-e^{-2\theta}\big)^{-\alpha/2}\big(1-e^{2(v-\theta)}\big)^{-1}{\rm d}\theta.
\end{align}
\end{theorem}
\noindent{\it Proof:} Since $Z$ is a symmetric $\alpha$-stable process, we have  for any $x\in \R^{d}$ and $b>0$
\[
\p_{x}\Big(b^{-1}Z_{\sigma_{b}^{+}}\in {\rm d} y;
\sigma_{b}^{+}<\infty\Big)=\p_{x/b}\Big(Z_{\sigma_{1}^{+}}\in {\rm d} y;
\sigma^{+}_{1}<\infty\Big),
\]
which implies that
\begin{equation}\label{scalint}
\p_{x}\Big(R_{\sigma_{e^{u}}^{+}}\in[e^{u},e^{u+\theta}];\sigma_{e^{u}}^{+}<\infty\Big)=\p_{e^{-u}x}\Big(R_{\sigma_{1}^{+}}\in[1,e^{\theta}];\sigma_{1}^{+}<\infty\Big).
\end{equation}
We first study the case $d=1$. Here, we assume that  $x=1$. From  (\ref{bl}), (\ref{scalint}) and Lemma 1, we have for $u, \theta\ge0$
\begin{align*}
P\Big(\xi_{T_{u}^{+}}\leq u+\theta;&\, T_{u}^{+}<\infty\Big)=\p_{e^{-u}}\left( R_{\sigma_{1}^{+}}\in[1,e^{\theta}];\sigma_{1}^{+}<\infty\right)\\
&=\frac1{\pi}\sin\left(\frac{\pi\alpha}{2}\right)(1-e^{-2u})^{\alpha/2}\int_{1\le |y|\le e^{\theta}}\big|1-|y|^{2}\big|^{-\alpha/2}|e^{-u}-y|^{-1}{\rm d} y,
\end{align*}
from which (\ref{teo4.1}) follows.\\
Now, we study the case $d\ge 2$. To this end, we fix $x\in\mathbb{R}^{d}$ such that $\|x\|=1$, and $w_{d}=2\pi^{d/2}\Big(\Gamma(d/2)\Big)^{-1}$. Hence using identity (\ref{bl}) and polar coordinates in $\R^d$, we have for $u, \theta\ge0$
\begin{align*}
\p_{e^{-u}x}&\left( R_{\sigma_{1}^{+}}\in[1,e^{\theta}];\sigma_{1}^{+}<\infty\right)\notag\\
&=\pi^{-(d/2+1)}\Gamma\left(\frac{d}{2}\right)\sin\left(\frac{\pi\alpha}{2}\right)(1-e^{-2u})^{\alpha/2}\int_{1\le\|y\|\le e^{\theta}}\big|1-\|y\|^{2}\big|^{-\alpha/2}\|e^{-u}x-y\|^{-d}{\rm d} y\notag\\
&=\pi^{-(d/2+1)}\Gamma\left(\frac{d}{2}\right)\sin\left(\frac{\pi\alpha}{2}\right)(1-e^{-2u})^{\alpha/2}\int_{1}^{e^{\theta}}{\rm d} r\frac{r^{d-1}}{(r^{2}-1)^{\alpha/2}}\notag\\
&\hspace{7.5cm}\times\int_{0}^{\pi}{\rm d}\phi\frac{w_{d-1}\sin^{d-2}\phi}{(r^{2}-2re^{-u}\cos\phi +e^{-2u})^{d/2}}.
\end{align*}
On the other hand, from formula 3.665 in \citep
{GR} we get for $r>1$
\[
\int_{0}^{\pi}{\rm d}\phi\frac{\sin^{d-2}\phi}{(r^{2}-2re^{-u}\cos\phi +e^{-2u})^{d/2}}=\frac{\pi^{1/2}\Gamma\big((d-1)/2\big)}{\Gamma(d/2)}e^{2u}r^{2-d}(r^2e^{2u}-1)^{-1},
\]
which implies that
\begin{align*}
\p_{e^{-u}x}&\left( R_{\sigma_{1}^{+}}\in[1,e^{\theta}];\sigma_{1}^{+}<\infty\right)\notag\\
&=\frac2{\pi}\sin\left(\frac{\pi\alpha}{2}\right)(1-e^{-2u})^{\alpha/2}e^{2u}\int_{1}^{e^{\theta}}{\rm d} r\,
r(r^{2}-1)^{-\alpha/2}(r^{2}-1)^{-1}.
\end{align*}
Therefore from Lemma 1 and (\ref{scalint}),  we conclude
\begin{align*}
P\Big(\xi_{T_{u}^{+}}&\leq u+\theta;
T_{u}^{+}<\infty\Big)\\
&=\frac2{\pi}\sin\left(\frac{\pi\alpha}{2}\right)(1-e^{-2u})^{\alpha/2}e^{2u}\int_{1}^{e^{\theta}}{\rm d} r\,
r(r^{2}-1)^{-\alpha/2}(r^{2}-1)^{-1},
\end{align*}
which proves (\ref{teo4.1}) for the case $d\ge 2$.

The second part of the theorem can be proved in a similar way. Indeed from the scaling property of $Z$, we have for $\theta\ge 0$ and $v\leq 0$
\begin{equation}\label{scalint1}
\p_{x}\Big(R_{\sigma_{e^{v}}^{-}}\in[e^{v-\theta},e^{v}];\sigma_{e^{v}}^{-}<\infty\Big)=\p_{e^{-v}x}\Big(R_{\sigma_{1}^{-}}\in[e^{-\theta},1];\sigma_{1}^{-}<\infty\Big).
\end{equation}
Assume that  $d=1$ and take  $x=1$. From  (\ref{bla}), (\ref{scalint1}) and Lemma 1, we have
\begin{align*}
P\Big(\xi_{T_{v}^{-}}\geq\theta- v;&\, T_{v}^{-}<\infty\Big)=\p_{e^{-v}}\left( R_{\sigma_{1}^{-}}\in[e^{-\theta}, 1];\sigma_{1}^{-}<\infty\right)\\
&=\frac1{\pi}\sin\left(\frac{\pi\alpha}{2}\right)(e^{-2v}-1)^{\alpha/2}\int_{e^{-\theta}\le|y|\le 1}\big|1-|y|^{2}\big|^{-\alpha/2}|e^{-v}-y|^{-1}{\rm d} y,
\end{align*}
from which (\ref{teo4.2}) follows.\\
Now, we study the case $d\ge 2$. To this end, we fix $x\in\mathbb{R}^{d}$ such that $\|x\|=1$, and set
$w_{d}=2\pi^{d/2}\Big(\Gamma(d/2)\Big)^{-1}$. Hence using (\ref{bla}), polar coordinates and formula 3.665 in \citep
{GR}, we get for $\theta\ge 0$ and $v\leq 0$
\begin{align*}
\p_{e^{-v}x}&\left( R_{\sigma_{1}^{-}}\in[e^{-\theta}, 1];\sigma_{1}^{-}<\infty\right)\notag\\
&=\pi^{-(d/2+1)}\Gamma\left(\frac{d}{2}\right)\sin\left(\frac{\pi\alpha}{2}\right)(e^{-2v}-1)^{\alpha/2}\int_{e^{-\theta}<\|y\|\le 1}\big|1-\|y\|^{2}\big|^{-\alpha/2}\|e^{-v}x-y\|^{-d}{\rm d} y\notag\\
&=\pi^{-(d/2+1)}\Gamma\left(\frac{d}{2}\right)\sin\left(\frac{\pi\alpha}{2}\right)(e^{-2v}-1)^{\alpha/2}\int_{e^{-\theta}}^1{\rm d} r\frac{r^{d-1}}{(1-r^{2})^{-\alpha/2}}\notag\\
&\hspace{7cm}\times\int_{0}^{\pi}{\rm d}\theta\frac{w_{d}\sin^{d-2}\theta}{(r^{2}+e^{-2v}-2re^{-v}\cos\theta)^{d/2}}\notag\\
&=\frac2{\pi}\sin\left(\frac{\pi\alpha}{2}\right)(e^{-2v}-1)^{\alpha/2}e^{-(2-d)v}\int_{e^{-\theta}}^1{\rm d} r\,
r^{d-1}(1-r^{2})^{-\alpha/2}(e^{-2v}-r^{2})^{-1}\notag
\end{align*}
Therefore from Lemma 1 and (\ref{scalint1}),  we conclude
\begin{align*}
P\Big(v-\xi_{T_{v}^{-}}&\leq \theta;
T_{u}^{-}<\infty\Big)\\
&=\frac{2}{\pi}\sin\left(\frac{\pi\alpha}{2}\right)(e^{-2v}-1)^{\alpha/2}e^{-(2-d)v}\int_{e^{-\theta}}^1{\rm d} r\
r^{d-1}(1-r^{2})^{-\alpha/2}(e^{-2v}-r^2)^{-1}.
\end{align*}
This complete the proof. \QED

Additional  computations yield the following corollary.
\begin{corollary}
Suppose that $\alpha<d$ and let $\underline{\xi}_{\infty}=\inf_{t\ge 0}\xi_t$. For $z\ge 0$,
\[
P\Big(-\underline{\xi}_{\infty}\in {\rm d} z\Big)=2\frac{\Gamma(d/2)}{\Gamma\big((d-\alpha)/2\big)\Gamma(\alpha/2)}e^{-(d-2)z}(e^{2z}-1)^{\alpha/2-1}{\rm d} z.
\]
\end{corollary}
\noindent{\it Proof:} We first note that
\[
\int_{0}^r u^{d-\alpha-1} (r^2-u^2)^{(\alpha-2)/2} {\rm d} u =\frac{r^{d-2}}{2}\frac{\Gamma(\alpha/2)\Gamma\big((d-\alpha)/2\big)}{\Gamma(d/2)},
\]
and that for $u\in[0,1]$ and $z>0$
\[
\int_0^{1-u^2}{\rm d} y y^{-\alpha/2}(e^{2z}-1+y)^{-1} (1-y-u^2)^{\alpha/2-1} =\frac{\pi}{\sin(\pi\alpha/2)}\frac{(e^{2z}-u^2)^{\alpha/2-1}}{(e^{2z}-1)^{\alpha/2}}.
\]
Thus, we have
\[
\begin{split}
\int_{0}^1&{\rm d} r\,
r^{d-1}(1-r^{2})^{-\alpha/2}(e^{2z}-r^2)^{-1}\\
&=\frac{2\Gamma(d/2)}{\Gamma(\alpha/2)\Gamma\big((d-\alpha)/2\big)}\int_{0}^1{\rm d} r\,
r(1-r^{2})^{-\alpha/2}(e^{2z}-r^2)^{-1}\int_{0}^r u^{d-\alpha-1} (r^2-u^2)^{(\alpha-2)/2} {\rm d} u\\
&=\frac{\Gamma(d/2)}{\Gamma(\alpha/2)\Gamma\big((d-\alpha)/2\big)}\int_{0}^1{\rm d} uu^{d-\alpha-1}\int_0^{1-u^2}{\rm d} y y^{-\alpha/2}(e^{2z}-1+y)^{-1} (1-y-u^2)^{\alpha/2-1} \\
&=\frac{\Gamma(d/2)}{\Gamma(\alpha/2)\Gamma\big((d-\alpha)/2\big)}\frac{\pi}{\sin(\pi\alpha/2)}(e^{2z}-1)^{-\alpha/2}
\int_{0}^1{\rm d} u\,u^{d-\alpha-1}(e^{2z}-u^2)^{\alpha/2-1}\\
&=\frac{\Gamma(d/2)}{\Gamma(\alpha/2)\Gamma\big((d-\alpha)/2\big)}
\frac{\pi}{2\sin(\pi\alpha/2)}(e^{2z}-1)^{-\alpha/2}e^{(d-2)z}
\int_{e^{2z}-1}^\infty{\rm d} r\,\frac{r^{\alpha/2-1}}{(r+1)^{d/2}}.
\end{split}
\]
Therefore, from the above computations and (\ref{teo4.2}) we get for $z> 0$
\[
\begin{split}
P\Big(\underline{\xi}_{\infty}\le-z\Big)&=P\big(T^-_{-z}<\infty\big)\\
&=\frac{2}{\pi}\sin\left(\frac{\pi\alpha}{2}\right)e^{-dz}\big(e^{2z}-1\big)^{\alpha/2}\int_{0}^{\infty}e^{-d\theta}\big(1-e^{-2\theta}\big)^{-\alpha/2}\big(1-e^{-2(z+\theta)}\big)^{-1}{\rm d}\theta\\
&=\frac{2}{\pi}\sin\left(\frac{\pi\alpha}{2}\right)e^{-(d-2)z}\big(e^{2z}-1\big)^{\alpha/2}\int_{0}^1{\rm d} r\,
r^{d-1}(1-r^{2})^{-\alpha/2}(e^{2z}-r^2)^{-1}\\
&=\frac{\Gamma(d/2)}{\Gamma(\alpha/2)\Gamma\big((d-\alpha)/2\big)}
\int_{e^{2z}-1}^\infty{\rm d} r\,\frac{r^{\alpha/2-1}}{(r+1)^{d/2}}.
\end{split}
\]
This complete the proof. \QED
\section{Entrance laws: points}
 For any $y\in \R$ and $r>0$, let
\[
T_{y}=\inf\{t>0:\xi_t=y\} \quad\textrm{ and }\quad
\sigma_{r}=\inf\{t>0: R_t=r\}.
\]
We also introduce
\[
{\rm P}^{\mu}_{\nu}(z)=\frac{1}{\Gamma(1-\mu)}\left(\frac{z+1}{z-1}\right)^{\mu/2}\, _{2}\mathcal{F}_{1}\left(-\nu, \nu+1;1-\mu;\frac{1-z}{2}\right)\quad z>1
\]
the so called Legendre function of the first kind.

The purpose of this section is to explicitly compute the probability that the process $\xi$ hits a point i.e.  $P(T_r<\infty)$, as well as some related quantities. Our study is based on the work of Port \citep
{po}, where the author computes the probability that the radial process $R$ hits a given point when $\alpha \in (1,2)$. We recall that  the radial process $R$  only hits points when $\alpha\in (1,2)$.

The one-point hitting probability for  $R$, presented in Port \citep
{po}  is given by the formula
\begin{equation}\label{port}
\p_x(\sigma_r<\infty)=\frac{2^{2-\alpha}\pi^{1/2}\Gamma\left((d+\alpha)/2-1\right)}{\Gamma\left((\alpha-1)/2\right)}r^{d/2+1-\alpha}\big|1-r^2\big|^{\alpha/2-1}{\rm P}^{1-d/2}_{-\alpha/2}\left(\frac{1+r^2}{|1-r^2|}\right),
\end{equation}
where $r>0$ and $x\in\R^d$ such that $\|x\|=1$. From  the Lamperti representation (\ref{lamp}) and  identity (\ref{port}), we obtain the one-point hitting problem  for $\xi$  as follows.
\begin{theorem}Let $1<\alpha<d$. Then for $y\in \R$
\[
P(T_y<\infty)=\frac{2^{2-\alpha}\pi^{1/2}\Gamma\left((d+\alpha)/2-1\right)}{\Gamma\left((\alpha-1)/2\right)}e^{(d/2-1)y}\big|e^{-2y}-1\big|^{\alpha/2-1}{\rm P}^{1-d/2}_{-\alpha/2}\left(\frac{1+e^{2y}}{|1-e^{2y}|}\right).
\]
\end{theorem}
\noindent{\it Proof:}  From the Lamperti representation (\ref{lamp}) of the  process $R$, we have for $y\in \R$ and $x\in \R^d$ satisfying $\|x\|=1$
\[
\p_x\big(\sigma_{e^y}<\infty\big)=P\left(\int_0^{T_y} e^{\alpha \xi_s} {\rm d} s <\infty\right).
\]
On the other hand, it is clear that
\begin{equation}\label{desigualdad}
T_y\exp\left\{\alpha\inf_{0\le u<T_y}\xi_u\right\}\le\int_0^{T_y} e^{\alpha \xi_s} {\rm d} s\le T_y\exp\left\{\alpha\sup_{0\le u<T_y}\xi_u\right\}.
\end{equation}
Hence if  $\int_0^{T_y} e^{\alpha \xi_s} {\rm d} s<\infty$ then we have that $T_y<\infty$, since the process $\xi$ drifts to $+\infty$ and $\inf_{0\le u<T_y}\xi_u>-\infty$.

Now, recall from Theorem 4 that the process $\xi$ does not creep upwards. If $T_{y}<\infty$, we have that the process $\xi$ makes a finite number of jumps across $y$ before time $T_{y}$ and then $\sup_{0\le u<T_y}\xi_u<\infty$. Hence from (\ref{desigualdad}), we deduce that  $\int_0^{T_y} e^{\alpha \xi_s} {\rm d} s<\infty$.  Therefore
\[
\p_x\big(\sigma_{e^y}<\infty\big)=P\left(T_y <\infty\right).
\]
This completes the proof. \QED
Now, we explore more elaborate hitting probabilities ($n$-point hitting problem) for the  L\'evy process $\xi$ when $1<\alpha<d$. This is possible  thanks to a result of Port \citep
{po} and  the  Lamperti representation (\ref{lamp}) of the process $R$. Let $B=\{r_{1},r_{2},\cdots,r_{n}\}$ where
$r_{1}<r_{2}<\cdots<r_{n}$.

Recall from \citep
{po}, that the potential density  $u(\cdot, \cdot)$ of the radial process $R$ which is specified by
\[
\e_z\left(\int_0^\infty\ind_{\{R_t\in A\}}{\rm d} t\right)=\frac{1}{2^{d/2}\Gamma(d/2+1)}\int_A {\rm d} y\, y^d u(\|z\|, y), \qquad\textrm{for}\quad z\in \R^d, A\in \mathcal{B}(\R_+),
\]
 satisfies (see Lemmas 2.1 and 2.2 in \citep
{po}), for $x,y>0$
\begin{equation}
 u(x,y) =
 \frac{2^{(d/2)-\alpha}\Gamma(d/2)\Gamma\big((d-\alpha)/2\big)}{\Gamma(\alpha/2)}(xy)^{1-d/2}|x^{2}-y^{2}|^{\alpha/2-1}{\rm P}_{-\alpha/2}^{1-d/2}\left(\frac{x^{2}+y^{2}}{|x^{2}-y^{2}|}\right)\notag,
\end{equation}
and
\begin{equation}
u(x,x)=\frac{\pi^{-1/2}2^{d/2-2}\Gamma((\alpha-1)/2)}{\Gamma(\big(\alpha+d\big)/2-1)}\frac{\Gamma(d/2)\Gamma\big((d-\alpha)/2\big)}{\Gamma(\alpha/2)}x^{\alpha-d},\notag
\end{equation}
and  that the matrix $U=\Big[u(r_{i},r_{j})\Big]_{n\times n}$ is invertible. Let us denote its inverse by $K_{B}=\Big[K_B(i,j)\Big]_{n\times n}$ and set   $\sigma_{B}=\inf\{t>0:R_{t}\in B\}$.

 According to Port, the probability that the process $R$ hits the set $B$ at a finite time is given by
\begin{equation}\label{eq1port}
\p_z(\sigma_{B}<\infty)=\sum_{i=1}^n \sum_{j=1}^n u(\|z\|, r_j)K_B(i,j),
\end{equation}
and the probability that it first hits the point $r_j$ is given by
\begin{equation}\label{eq2port}
\mathbb{P}_{z}\Big(R_{\sigma_{B}} =r_{j}; \, \sigma_{B}
< \infty\Big)=\sum_{i=1}^{n}u(\|z\|,r_{i})K_{B}(i,j).
\end{equation}
For a two point set $B=\{r_{1},r_{2}\}$ we have that
\begin{align*}
K_{B}=\frac{1}{\Delta}\left(
\begin{array}{cc}
U_{22} & -U_{12}\\
-U_{12} & U_{11}
 \end{array}
\right),
\end{align*}
where $\Delta=U_{11}U_{22}-U_{12}^{2}$. Then from (\ref{eq1port}) and (\ref{eq2port}), we have
\[
\p_z(\sigma_{B}<\infty)=\frac{u(\|z\|,r_{1})u(r_{2},r_{2})+u(\|z\|,r_{2})u(r_{1},r_{1})}{u(r_{1},r_{1})u(r_{2},r_{2})-u(r_{1},r_{2})^{2}}-\frac{u(r_{1}, r_2)[u(\|z\|,r_{1})+u(\|z\|,r_{2})]}{u(r_{1},r_{1})u(r_{2},r_{2})-u(r_{1},r_{2})^{2}},
\]
and
\begin{align*}
\p_z(\sigma_{r_1}<\sigma_{r_2})&=\frac{u(\|z\|,r_{1})u(r_{2},r_{2})-u(\|z\|,r_{2})u(r_{2},r_{1})}{u(r_{1},r_{1})u(r_{2},r_{2})-u(r_{1},r_{2})^{2}},\notag\\
\p_z(\sigma_{r_2}<\sigma_{r_1})&=\frac{u(\|z\|,r_{2})u(r_{1},r_{1})-u(\|z\|,r_{1})u(r_{1},r_{2})}{u(r_{1},r_{1})u(r_{2},r_{2})-u(r_{1},r_{2})^{2}}.
\end{align*}

Hence the two-point hitting probabilities for the L\'evy process $\xi$ are as follows.
\begin{theorem}
Suppose that $1<\alpha<d$ and fix $-\infty<v<0<u<\infty$. Define
\[
T_{\{v,u\}}=\inf\{t>0 : \xi_t\in\{v,u\}\}.
\]
We have
\[
P\Big(T_{\{v,u\}}<\infty\Big)=\frac{u(1,e^v)u(e^u,e^u)+u(1,e^u)u(e^v,e^v)}{u(e^v,e^v)u(e^u,e^u)-u(e^v,e^u)^{2}}-\frac{u(e^v, e^u)[u(1,e^v)+u(1,e^u)]}{u(e^v,e^v)u(e^u,e^u)-u(e^v,e^u)^{2}},
\]
\[
 P\Big(\xi_{   T_{\{v,u\}}  } =v\Big) = f(1, e^v, e^u)\qquad\textrm{
and}\qquad
 P\Big(\xi_{T_{\{v,u\}}} =u\Big) = f(1, e^u, e^v),
 \]
where
\[
f(x,a,b) = \frac{\frac{u(x,a)}{u(b,a)} - \frac{u(x,b)}{u(b,b)}}{\frac{u(a,a)}{u(b,a)} - \frac{u(a,b)}{u(b,b)}}.
\]
\end{theorem}

\section{Wiener-Hopf factorization.}
In this section we work in the polar case and  compute explicitly the characteristic exponent of the process $\xi$ using its Wiener-Hopf factorization. Denote by $\{(L_{t}^{-1},H_{t}):t\geq 0\}$ and $\{(\widehat{L}_{t}^{-1},%
\widehat{H}_{t}):t\geq 0\}$ the (possibly killed) bivariate subordinators
representing the ascending and descending ladder processes of $\xi$ (see  \citep
{Be} for a proper definition). Write $%
\kappa (\theta,\lambda )$ and $\widehat{\kappa }(\theta,\lambda)$ for their joint
Laplace exponents for $\theta,\lambda\geq 0$. For convenience we will write
\[
\widehat{\kappa} (0,\lambda
 )=\widehat{q}+\widehat{\mathrm{c}}\lambda
 +\int_{(0,\infty
)}(1-e^{-\lambda
 x})\Pi _{\widehat{H}}({{\rm d}}x),
\]%
where $\widehat{q}\geq 0$ is the killing rate of $\widehat{H}$ so that $\widehat{q}>0$ if and only if $%
\lim_{t\uparrow \infty }\xi_{t}=\infty $, $\widehat{\mathrm{c}}\geq 0$ is the drift of $%
\widehat{H}$ and $\Pi _{\widehat{H}}$ is its jump measure.  Similar notation will also be used for $\kappa (0,\lambda
 )$ by replacing  $\widehat{q}$, $\widehat{\xi }$, $\widehat{\mathrm{c}}$ and $%
\Pi _{\widehat{H}}$ by $q$, $\xi $, $\mathrm{c}$ and $%
\Pi _{H}$. Note that necessarily $q=0$ since $\lim_{t\uparrow \infty }\xi_{t}=\infty$.

Associated with the ascending and descending ladder processes are the
bivariate renewal functions $V$ and $\widehat{V}$. The
former is defined by
\[
V({{\rm d}}s, {{\rm d}}x)=\int_{0}^{\infty }{\rm d}t\cdot
P(L_{t}^{-1}\in {{\rm d}}s, H_{t}\in {{\rm d}}x)
\]%
and taking double Laplace transforms shows that
\begin{equation}
\int_{0}^{\infty }\int_{0}^{\infty }e^{-\theta s-\lambda
 x}V({
{\rm d}}s,{{\rm d}}x)=\frac{1}{\kappa (\theta ,\lambda )}\quad\text{for }\theta ,\lambda
\geq 0  \label{doubleLT}
\end{equation}%
with a similar definition and relation holding for $\widehat{V}$.
These bivariate renewal measures are essentially the Green's measures of
the ascending and descending ladder processes. With an abuse of notation we shall also write $V({{\rm d}}x)$ and $%
\widehat{V}({{\rm d}}x)$ for the marginal measures $V(%
[0,\infty ), {{\rm d}}x)$ and $\widehat{V}([0,\infty ), {{\rm d}}x)
$ respectively. (Since we shall never use the marginals  $V(%
{{\rm d}}s,[0,\infty ))$ and $\widehat{V}({{\rm d}}s,[0,\infty ))
$ there should be no confusion). Note that local time at the maximum is defined only up to a
multiplicative constant. For this reason, the exponent $\kappa $ can only be
defined up to a multiplicative constant and hence the same is true of the
measure $V$ (and then obviously this argument applies to $\widehat{%
V}$).

The main result of this section is the Wiener-Hopf factorization of the characteristic exponent of the L\'evy process $\xi$.
\begin{theorem}
Let $\alpha<d$ and $\xi$ be the hypergeometric-stable L\'evy process. Then its characteristic exponent $\Psi$ enjoys the following Wiener-Hopf factorization
\begin{equation}\label{whexp}
\begin{split}
\Psi(\lambda)&=2^{\alpha}\frac{\Gamma((-i\lambda+\alpha)/2)}{\Gamma(-i\lambda/2)}\frac{\Gamma((i\lambda +d)/2)}{\Gamma((i\lambda +d-\alpha)/2)}\\
&=2^{\alpha}\frac{\Gamma(d/2)\Gamma((-i\lambda+\alpha)/2)}{\Gamma((d-\alpha)/2)\Gamma(-i\lambda/2)}\times\frac{\Gamma((d-\alpha)/2)\Gamma((i\lambda +d)/2)}{\Gamma(d/2)\Gamma((i\lambda +d-\alpha)/2)}
\end{split}
\end{equation}
where the first equality hold up to a multiplicative constant.
\end{theorem}

The proof of Theorem 7 relies on the  computation of the Laplace exponents of the ascending ladder height  and the descending  ladder height processes of $\xi$.
\begin{lemma}
Let $\alpha<d$ and $\xi$ be the hypergeometric-stable L\'evy process. The Laplace exponent of its descending ladder height process
$\widehat{H}$  is given by
\begin{equation}\label{dwh}
\hat{\kappa}(0,\lambda)=\frac{\Gamma((d+\lambda)/2)\Gamma((d-\alpha)/2)}{\Gamma(d/2)\Gamma((d-\alpha+\lambda)/2)}.
\end{equation}
\end{lemma}
\noindent {\it Proof:} Recall from the proof of Corollary 2 that
\[
P\left(-\inf_{t\geq0}\xi_{t}\le z\right)=\frac{\Gamma(d/2)}{\Gamma((d-\alpha)/2)\Gamma(\alpha/2)}\int_{0}^{e^{2z}-1}(u+1)^{-d/2}u^{\alpha/2-1}{\rm d} u.
\]
Also recall that $\widehat{V}$ denotes the renewal  function  associated with $\widehat{H}$. From  Proposition VI.17 in \citep
{Be}, we know that
\[
\widehat{V}(z):=\widehat{V}([0,z])=\widehat{V}([0,\infty))P\left(-\inf_{t\geq0}\xi_{t}\leq
z\right) \qquad\textrm{for all }\quad z\ge 0.
\]
As we mentioned before, it is  well known that $\widehat{V}$ is unique up to a multiplicative constant which depends on the normalization of local time of $\xi$ at its infimum. Without loss of generality we may therefore assume in the forthcoming analysis that $\widehat{V}(\infty)$, which is equal to the reciprocal of killing rate of the descending ladder height process, may be taken identically equal to 1. Hence
\[
\widehat{V}(z)=\frac{\Gamma(d/2)}{\Gamma((d-\alpha)/2)\Gamma(\alpha/2)}\int_{0}^{e^{2z}-1}(u+1)^{-d/2}u^{\alpha/2-1}{\rm d} u.
\]
Now, let $K(\alpha,d)=\Gamma(d/2)\big(\Gamma((d-\alpha)/2)\Gamma(\alpha/2)\big)^{-1}$ and note
\begin{align*}
\lambda\int_{0}^{\infty}e^{-\lambda x}\widehat{V}(x){\rm d} x&=\lambda
K(\alpha,d)\int_{0}^{\infty}{\rm d} x \,e^{-\lambda
x}\int_{0}^{e^{2x}-1}{\rm d} u\, (u+1)^{-d/2}u^{\alpha/2-1}\notag\\
&=K(\alpha,d)\int_{0}^{\infty}(u+1)^{-(d+\lambda)/2}u^{\alpha/2-1}{\rm d} u\notag\\
&=K(\alpha,d)\int_0^\infty u^{(d-\alpha+\lambda)/2-1}(1-u)^{\alpha/2-1} {\rm d} u\\
&=\frac{\Gamma(d/2)\Gamma((d+\lambda-\alpha)/2)}{\Gamma((d+\lambda)/2)\Gamma((d-\alpha)/2)}.
\end{align*}
Finally, from (\ref{doubleLT}) we deduce that
\[
\hat{\kappa}(0,\lambda)=\frac{\Gamma((d+\lambda)/2)\Gamma((d-\alpha)/2)}{\Gamma(d/2)\Gamma((d-\alpha+\lambda)/2)}.
\]
This completes the proof.\QED
For the computation of  the Laplace exponent of the ascending ladder height process
$H$, we will make  use of an important identity obtained by
Vigon \citep
{vi} that relates $\Pi_H$, the L\'evy measure of the ascending ladder height process $H$,
with that of the L\'evy process $\xi$ and $\widehat{V}$, the potential measure of the descending ladder height process $\widehat{H}$. Specifically, defining $\overline{\Pi}_H(x)=\Pi_{H}(x,\infty)$, the identity states that
\begin{equation}\label{vi}
\overline{\Pi}_H(r)=\int_{0}^{\infty}\widehat{V}({\rm d} l)\overline{\Pi}^{+}(l+r) \qquad r>0,
\end{equation}
where $\overline{\Pi}^+(u)=\Pi(u,\infty)$ for $u>0$.

Now,  recall the following property of the hypergeometric
function $_2\mathcal{F}_{1}$ (see for instance identity (3.1.9) in \citep
{aak})
\begin{equation}\label{p2}
_2\mathcal{F}_{1}(a,b;a-b+1;x)=(1+x)^{-a}\, _2\mathcal{F}_{1}\left(a/2,(a+1)/2;a-b+1;\frac{4x}{(1+x)^{2}}\right),
\end{equation}
and note that  the L\'evy measure of the process $\xi$ can be written as follows
\[
\begin{split}
\Pi({\rm d} y)&=\frac{e^{-\alpha
y}}{(1+e^{-2y})^{(\alpha+d)/2}}\overline{F}\left(\frac{4e^{-2y}}{(1+e^{-2y})^{2}}\right)\ind_{\{y>0\}}{\rm d} y\\
& \hspace{4cm}+\frac{e^{d
y}}{(1+e^{2y})^{\alpha+d/2}}\overline{F}\left(\frac{4e^{2y}}{(1+e^{2y})^{2}}\right)\ind_{\{y<0\}}{\rm d} y.
\end{split}
\]
Therefore
\begin{equation}\label{c3}
\begin{split}
\Pi({\rm d} y)&=\frac{2^{\alpha}\alpha(d/2)_{\alpha/2}}{\Gamma(1-\alpha/2)}e^{-\alpha
y}\, _2\mathcal{F}_{1}\Big((\alpha+d)/2,\alpha/2+1;d/2;e^{-2y}\Big)\ind_{\{y>0\}}{\rm d} y\\
&\hspace{3cm}+\frac{2^{\alpha}\alpha(d/2)_{\alpha/2}}{\Gamma(1-\alpha/2)}e^{d y}\, _2\mathcal{F}_{1}\Big(\alpha+d/2,\alpha/2+1;d/2;e^{2y}\Big)\ind_{\{y<0\}}{\rm d} y.
\end{split}
\end{equation}
\begin{lemma}
Let $\alpha<d$ and $\xi$ be the hypergeometric-stable L\'evy process.
The Laplace exponent of its ascending ladder height process
$H$  is given by
\begin{equation}\label{awh}
\kappa(0,\lambda)=\frac{2^{\alpha}\Gamma(d/2)\Gamma((\lambda+\alpha)/2)}{\Gamma((d-\alpha)/2)\Gamma(\lambda/2)}.
\end{equation}
\end{lemma}
\noindent{\it Proof:} We first note from the proof of Lemma 2, that  the
renewal measure $\widehat{V}({\rm d} y)$ associated with $\widehat{H}$ satisfies
\begin{equation}\label{renmdown}
\widehat{V}({\rm d} y)=\frac{2\Gamma(d/2)}{\Gamma((d-\alpha)/2)\Gamma(\alpha/2)}e^{(2-d)y}(e^{2y}-1)^{\alpha/2-1}{\rm d} y.
\end{equation}
We also recall the following property of the Gamma function,
 $$\Gamma(1-\alpha/2)\Gamma(\alpha/2)=\frac{\pi}{\sin(\pi \alpha/2)}.$$
From Vigon's formula (\ref{vi}) and identity (\ref{c3}), we have
\[
\begin{split}
\overline{\Pi}_H(x)&=\frac{2^{\alpha+1}\alpha \sin (\alpha\pi/2)}{\pi}\frac{\Gamma((d+\alpha)/2)}{\Gamma((d-\alpha)/2)}\int_{0}^{\infty}{\rm d} y\, e^{(2-d)y}(e^{2y}-1)^{\alpha/2-1}\\
&\hspace{5.5cm}\times\int_{x+y}^{\infty}{\rm d} u\, e^{-\alpha
u}\, _2\mathcal{F}_{1}\Big((\alpha+d)/2,\alpha/2+1;d/2;e^{-2u}\Big).
\end{split}
\]
On the other hand from the definition of $_2\mathcal{F}_1$, we get
\[
\begin{split}
\int_{x+y}^{\infty}{\rm d} u\, e^{-\alpha
u}\, _2\mathcal{F}_{1}&\Big((\alpha+d)/2,\alpha/2+1;d/2;e^{-2u}\Big)\\
&=\frac{1}{2}\int_{0}^{e^{-2(x+y)}}{\rm d} z\,z^{\alpha/2-1}\,_2\mathcal{F}_{1}\Big((\alpha+d)/2,\alpha/2+1;d/2;z\Big)\notag\\
&=\frac{e^{-\alpha(x+y)}}{\alpha}\, _2\mathcal{F}_1\Big((d+\alpha)/2,\alpha/2;d/2;e^{-2(x+y)}\Big).
\end{split}
\]
Set
\[
C(\alpha, d)=\frac{2^{\alpha+1} \sin (\alpha\pi/2)}{\pi}\frac{\Gamma((d+\alpha)/2)}{\Gamma((d-\alpha)/2)}.
\]
Hence putting the pieces together,  we obtain
\begin{align}
\overline{\Pi}_H(x)&=C(\alpha, d)e^{-\alpha x}\int_{0}^{\infty}\, _2\mathcal{F}_{1}\big((d+\alpha)/2,\alpha/2;d/2;e^{-2(x+y)}\big)e^{y(2-d-\alpha)}(e^{2y}-1)^{\alpha/2-1}{\rm d} y\notag\\
&=C(\alpha, d)\sum_{k=0}^{\infty}e^{-2x(\alpha/2+k)}\frac{((d+\alpha)/2)_{k}(\alpha/2)_{k}}{(d/2)_{k}k!}\int_{0}^{\infty}e^{-2y(d/2+k)}(1-e^{-2y})^{\alpha/2-1}{\rm d} y\notag\\
&=\frac{C(\alpha, d)}2\sum_{k=0}^{\infty}e^{-2x(\alpha/2+k)}\frac{((d+\alpha)/2)_{k}(\alpha/2)_{k}}{(d/2)_{k}k!}\int_{0}^{1}u^{d/2+k-1}(1-u)^{\alpha/2-1}{\rm d} u\notag\\
&=\frac{C(\alpha, d)}2\sum_{k=0}^{\infty}e^{-2x(\alpha/2+k)}\frac{((d+\alpha)/2)_{k}(\alpha/2)_{k}}{(d/2)_{k}k!}\frac{\Gamma(d/2+k)\Gamma(\alpha/2)}{\Gamma((d+\alpha)/2+k)}\notag\\
&=\frac{C(\alpha, d)}2\frac{\Gamma(d/2)\Gamma(\alpha/2)}{\Gamma((d+\alpha)/2)}e^{-\alpha
x}\sum_{k=0}^{\infty}e^{-2k x}\frac{(\alpha/2)_{k}}{k!}\notag\\
&=\frac{2^{\alpha} \sin (\alpha\pi/2)}{\pi}\frac{\Gamma(d/2)\Gamma(\alpha/2)}{\Gamma((d-\alpha)/2)}e^{-\alpha
x}(1-e^{-2x})^{-\alpha/2}\notag.
\end{align}
From Theorem 3, we deduce  that the process $\xi$ does not creep upwards. Hence by Theorem VI.19 of \citep
{Be} the ascending ladder height process $H$ has no drift. Also recall that the process $\xi$ drift to $\infty$ which implies that the process  $H$ has no killing term. Therefore the Laplace exponent $\kappa(0,\lambda)$ of $H$ is given by
\begin{align*}
\frac{\kappa(0,\lambda)}{\lambda}=\frac{2^{\alpha} \sin (\alpha\pi/2)}{\pi}\frac{\Gamma(d/2)\Gamma(\alpha/2)}{\Gamma((d-\alpha)/2)}\int_{0}^{\infty}e^{-\lambda
x}e^{-\alpha
x}(1-e^{-2x})^{-\alpha/2}{\rm d} x.
\end{align*}
By integrating by parts and a change of variable, we get
\begin{align}
\kappa(0,\lambda)&=\frac{\alpha 2^{\alpha} \sin (\alpha\pi/2)}{\pi}\frac{\Gamma(d/2)\Gamma(\alpha/2)}{\Gamma((d-\alpha)/2)}\int_{0}^{\infty}\Big(1-e^{-(\lambda/2)
x}\Big)\frac{e^{x}}{(e^{x}-1)^{\alpha/2+1}}{\rm d} x.\notag
\end{align}
According to Theorem 3.1 of \citep
{cpp} (see Theorem 3.1), the previous integral satisfies
\[
\int_{0}^{\infty}\Big(1-e^{-(\lambda /2)x
}\Big)\frac{e^{x}}{(e^{x}-1)^{\alpha/2+1}}{\rm d} x=-\frac{\Gamma(-\alpha/2)\Gamma((\lambda+\alpha)/2)}{\Gamma(\lambda/2)},
\]
where $\Gamma(-\alpha/2)=-\alpha^{-1}\Gamma(1-\alpha/2)$. Therefore,
\[
\kappa(0,\lambda)=\frac{2^{\alpha}\Gamma(d/2)\Gamma((\lambda+\alpha)/2)}{\Gamma((d-\alpha)/2)\Gamma(\lambda/2)}
\]
This completes the proof.\QED
\noindent{\it Proof of Theorem 7:}
From the fluctuation theory of L\'evy processes, it is known that  Wiener-Hopf factorization of the characteristic exponent of $\xi$ is given by
\[
\psi(\lambda)=\kappa(0,-i\lambda)\times \hat{\kappa}(0,i\lambda)
\]
up to a multiplicative constant. Hence, the result follows  from Lemmas 2 and 3.
\QED
\begin{rem}

  We have obtained the characteristic exponent for the process $\xi$ in the case where
   $\alpha<d$ using the Wiener-Hopf factorization.
   We will now see that the same formula holds true in the example studied in \citep
{CPY}: $\alpha=d=1$.

    Recall that they obtained the following characteristic exponent  of  $\xi$:
     $$E \Big[\exp \{i\lambda\xi_t\}\Big] = \exp \left\{-t\lambda \tanh\left(\frac{\pi\lambda}{2}\right)\right\}, \quad t\geq 0, \quad \lambda \in
     R.$$
We have 
$$\psi(\lambda)=\lambda \tanh\left(\frac{\pi\lambda}{2}\right)=\frac{\frac{\pi}{\cosh(\pi \lambda/2)}}{\frac{\pi}{ (\lambda/2) \sinh(\pi \lambda/2)}}=\frac{\vert \Gamma\left( \frac{i\lambda+1}{2}\right) \vert ^2}{\vert \Gamma\left( \frac{i\lambda}{2}\right) \vert^2}
=\left(\frac{i\lambda+1}{2}\right)_{1/2}\left(-\frac{i\lambda}{2}\right)_{1/2}.$$
Recall that the characteristic exponent in the case $\alpha<d$ is given by (\ref{whexp}).  From the above computation we note  that this formula still holds for the case $\alpha=d=1$.

From the unicity of the Wiener-Hopf factorization, we  deduce that  the characteristic exponent of the subordinators $\hat{H}$ and $H$ are:
 \[
 \hat{\kappa}(0,i\lambda)=\left(\frac{i\lambda+1}{2}\right)_{1/2}\qquad\qquad\kappa(0,-i\lambda)=\left(-\frac{i\lambda}{2}\right)_{1/2}.
 \]

 \end{rem}

\section{$n$-tuple laws at first and last passage times.}

Recall that the renewal measure $\widehat{V}(\ud y)$ associated with $\widehat{H}$ satisfies
\[
\widehat{V}(\ud y)=\frac{2\Gamma(d/2)}{\Gamma((d-\alpha)/2)\Gamma(\alpha/2)}e^{(2-d)y}(e^{2y}-1)^{\alpha/2-1}\ud y.
\]
From the form of the Laplace exponent of $H$ and (\ref{doubleLT}), we get that the renewal measure $V(\ud y)$ associated with $H$ satisfies
\[
V(\ud y)=\frac{\Gamma((d-\alpha)/2)}{2^{\alpha-1}\Gamma(d/2)\Gamma(\alpha/2)}(1-e^{-2y})^{\alpha/2-1}\ud y.
\]
Since we have explicit expressions for the renewal functions $V$ and $\widehat{V}$,  we can get, from the main results of Doney and Kyprianou \citep
{DK} and Kyprianou et al. \citep
{kpr},   $n$-tuple  laws at first and last passage times  for  the L\'evy process $\xi$ and the radial part of the symmetric stable L\'evy process $Z$.

Marginalizing the quintuple law at first passage of Doney and Kyprianou \citep
{DK} (see Theorem 3) and by the  Lamperti representation (\ref{lamp}), we now obtain the following new identities.

\begin{proposition} Let $\overline{\xi}_t=\sup_{0\le s\le t} \xi_s$. For $y\in [0,x]$, $v\ge y$ and $u>0$,
\[
\begin{split}
P&\Big(\xi_{T^+_x}-x\in {\rm d} u, x-\xi_{T^+_x-}\in {\rm d} v, x-\overline{\xi}_{T^+_x-}\in {\rm d} y\Big)\\
&=\frac{4\alpha\Gamma((\alpha+d)/2)}{\Gamma(d/2)\Gamma(\alpha/2)}\frac{\sin(\alpha\pi/2)}{\pi}(1-e^{-2(x-y)})^{\alpha/2-1}e^{(2-d)(v-y)}(e^{2(v-y)}-1)^{\alpha/2-1}\\
&\hspace{5cm}\times e^{-\alpha(u+v)} \ _{2}\mathcal{F}_{1}\Big((\alpha+d)/2), \alpha/2+1; d/2;e^{-2(u+v)}\Big){\rm d} y {\rm d} v {\rm d} u.
\end{split}
\]
For $z\in [x,1]$, $w\in [0,z]$ and $\theta>1$
\[
\begin{split}
\p_x&\left(\sup_{0\le s <\sigma^+_1}R_s\in {\rm d} z , R_{\sigma^+_1-}\in {\rm d} w, R_{\sigma^+_1}\in {\rm d} \theta\right)\\
&=\frac{4\alpha\Gamma((\alpha+d)/2)}{\Gamma(d/2)\Gamma(\alpha/2)}\frac{\sin(\alpha\pi/2)}{\pi}z^{3-d-\alpha}w^{d-1}\theta^{-\alpha-2}(z^2-x^2)^{\alpha/2-1}\\
&\hspace{4cm}\times (z^2-w^2)^{\alpha/2-1} \ _{2}\mathcal{F}_{1}\Big((\alpha+d)/2), \alpha/2+1; d/2;(w/\theta)^2\Big){\rm d} z {\rm d} w {\rm d} \theta.
\end{split}
\]
\end{proposition}
Note that the normalizing constant above is chosen to make the densities on the right-hand
side  distributions. It is also important to remark that the triple law for the L\'evy process $\xi$ extends  the identity in (\ref{teo4.1}).

Let us define the last passage time and the future infimum for the processes $\xi$ and $R$, respectively
\[
U_x=\sup\{t: \xi_t<x \},\quad L_x=\sup\{t: R_t<x \},\quad J_t=\inf_{s\ge t}\xi_s\quad \textrm{ and }\quad  F_t=\inf_{s\ge t}R_s.
\]
From Proposition 2.3 in Millar \citep
{mi}, we know that if $z>0$ the radial process $R$ of the symmetric stable L\'evy process is regular for both $(z,\infty)$ and $[0,z)$. Hence from the Lamperti representation (\ref{lamp}), we  deduce that the L\'evy process $\xi$ is regular for both $(-\infty, 0)$ and $(0,\infty)$. Now, applying corollaries 2 and 5 in Kyprianou et al. \citep
{kpr}, we obtain quadruple laws at last passage times for  $\xi$ and $R$.
\begin{proposition} For $x,v>0$, $0\leq y<x+v$ and $w\ge v>0$,
\[
\begin{split}
P&\Big(-J_0\in {\rm d} v, J_{U_x}-x\in {\rm d} u, x-\xi_{U_x-}\in {\rm d} y, \xi_{U_x}-x\in {\rm d} w\Big)\\
&=\frac{8\alpha\Gamma((\alpha+d)/2)}{\Gamma((d-\alpha)/2)\Gamma^2(\alpha/2)}\frac{\sin(\alpha\pi/2)}{\pi}e^{(2-d)(v+w-u)}(e^{2v}-1)^{\alpha/2-1}(e^{2(w-u)}-1)^{\alpha/2-1}\\
&\hspace{1cm}\times (1-e^{-2(x+v-y)})^{\alpha/2-1}e^{-\alpha(w+y)}\ _{2}\mathcal{F}_{1}\Big((\alpha+d)/2), \alpha/2+1; d/2;e^{-2(w+y)}\Big){\rm d} w {\rm d} y {\rm d} u {\rm d} v.
\end{split}
\]
For $x,b>0$, we have on $v\ge x^{-1}\lor b^{-1}$, $v^{-1}<y<b$ and $b<u\le w<\infty$
\[
\begin{split}
\p_x&\Big(1/F_0\in {\rm d} v, R_{L_b-}\in {\rm d} y, R_{L_b}\in {\rm d} w, F_{L_b}\in {\rm d} u\Big)\\
&=\frac{8\alpha\Gamma((\alpha+d)/2)}{\Gamma((d-\alpha)/2)\Gamma^2(\alpha/2)}\frac{\sin(\alpha\pi/2)}{\pi}b^{d-2\alpha}v^{1-d}yw^{1-d-\alpha}u^{d-\alpha-1}(v^2-1)\big(y^2-(bv)^{-2}\big)^{\alpha/2-1}\\
&\hspace{2.8cm}\times \big(w^2-(bu)^{2}\big)^{\alpha/2-1} \ _{2}\mathcal{F}_{1}\Big((\alpha+d)/2), \alpha/2+1; d/2;(y/bw)^{2}\Big){\rm d} v {\rm d} y {\rm d} w {\rm d} u.
\end{split}
\]
\end{proposition}
We conclude this section with a nice formula for the potential kernel of the L\'evy process $\xi$ killed as it enters $(-\infty,0)$, that follows from Theorem VI.20 in Bertoin \citep
{Be}.
\begin{proposition} There exist a constant $k>0$ such that for every  measurable function $f:[0,\infty)\to [0,\infty)$ and $x\ge 0$, one has
\[
\begin{split}
E_x&\left(\int_0^{T^-_{0}}f(\xi_t){\rm d} t\right)\\
&\hspace{1cm}=k\frac{2^{2-\alpha}}{\Gamma^2(\alpha/2)}\int_0^\infty {\rm d} y (1-e^{-2y})^{\alpha/2-1} \int_{0}^x {\rm d} z e^{(2-d)z}(e^{2z}-1)^{\alpha/2-1}f(x+y-z).
\end{split}
\]
In particular, the potential measure of the L\'evy process $\xi$ killed as it enters $(-\infty,0)$ has a density which is given by
\[
r(x,u)=k\frac{2^{2-\alpha}}{\Gamma^2(\alpha/2)}\int_{(u-x)\lor 0}^{u}   (1-e^{-2y})^{\alpha/2-1}e^{(2-d)(x+y-u)}(e^{2(x+y-u)}-1)^{\alpha/2-1} {\rm d} y.
\]
\end{proposition}
Note that from the previous proposition, we can obtain the potential kernel of the radial process $R$ killed as it enters $(0,1)$. Let $x>1$, then
 \[
\begin{split}
\e_x&\left(\int_0^{\sigma^-_{1}}f(R_t){\rm d} t\right)=E_{\log x}\left(\int_0^{T^-_{0}}f(e^{\xi_t})e^{\alpha \xi_t}{\rm d} t\right)\\
&=k\frac{2^{2-\alpha}}{\Gamma^2(\alpha/2)}\int_0^\infty {\rm d} y (1-e^{-2y})^{\alpha/2-1} \int_{0}^{\log x} {\rm d} z e^{(2-d)z}(e^{2z}-1)^{\alpha/2-1}x^{\alpha}e^{\alpha(y-z)}f(xe^{y-z}).
\end{split}
\]
In particular,
 \[
\begin{split}
\e_x&\left(\sigma^-_{1}\right)=E_{\log x}\left(\int_0^{T^-_{0}}e^{\alpha \xi_t}{\rm d} t\right)\\
&\hspace{.5cm}=k\frac{2^{2-\alpha}}{\Gamma^2(\alpha/2)}\int_0^\infty {\rm d} y (1-e^{-2y})^{\alpha/2-1} \int_{0}^{\log x} {\rm d} z e^{(2-d)z}(e^{2z}-1)^{\alpha/2-1}x^{\alpha}e^{\alpha(y-z)}\\
&\hspace{.5cm}=k\frac{x^{\alpha}}{2\Gamma(\alpha)} \int_{x^{-2}}^{1} {\rm d} u \,u^{d/2-1}(1-u)^{\alpha/2-1}.
\end{split}
\]

\noindent {\bf Acknowledgements.} This research was supported by EPSRC grant EP/D045460/1, CONACYT
grant (200419), and the project PAPIITT-IN120605. We are much
indebted to Andreas Kyprianou for many fruitful
discussions on L\'evy processes and fluctuation theory and to Marc Yor for pointing out the relationship with  their pioneering work \citep
{CPY} as well as for many enlighting related conversations.



\begin{thebibliography}{99}
\bibitem{aak} \sc Andrews, G.E.,  Askey, R., and Roy, R. (1999). {\it Special
Functions.} \rm Cambridge University Press, Cambridge.

\bibitem{Be} \sc Bertoin, J. (1996). {\it L\'evy Processes.} \rm
Cambridge University Press, Cambridge.


\bibitem{Bl} \sc Blumenthal, R., Getoor, R. K., and Ray, D. B. (1961).
\rm On the distribution of first hits for the symmetric stable
processes. {\it Trans. Amer. Math. Soc.},  {\bf 99},  540--554.



\bibitem{cpp} \sc Caballero, M. E., Pardo, J. C., and P\'erez, J. L. (2009). \rm On Lamperti Stable Processes. {\it To appear in Probability and  Mathematical  Statisistics}.

\bibitem{CPY} \sc Carmona, P., Petit, F., and Yor, M. (2001). \rm
Exponential Functionals of L\'evy processes. {\it L\'evy Processes, Theory and Applications, Eds. O.E.~Barndorff Nielsen et al.} Birkhauser , 41--56.


\bibitem{ckp} \sc Chaumont, L., Kyprianou, A. E., and Pardo J. C. (2009). \rm Some explicit identities
associated with positive self-similar Markov processes. \textit{ Stoch. Process.  Appl.,} {\bf 119}, 980-1000. 

\bibitem{Do} \sc Doney, R. A. (2007). \it Fluctuation theory for L\'evy processes.
\rm Ecole d'\'et\'e de Probabilit\'es de Saint-Flour, Lecture Notes
in Mathematics No. 1897. \rm Springer.

\bibitem{DK}\sc Doney, R.A. and Kyprianou, A.E. (2006). \rm Overshoots and undershoots of L\'evy processes.  {\it Ann.  Appl. Probab., } 16(1), 91--106.



\bibitem{GR}\sc Gradshtein, I.S. and Ryshik, I.M. (2007). {\it Table of Integrals, Series and Products },\rm Academic Press, San Diego.

\bibitem{kpr}\sc  Kyprianou, A.E.,  Pardo, J.C. and  Rivero, V. (2009). \rm Exact and asymptotic n-tuple laws at first and last passage.  {\it To appear in Ann. Appl. Probab. }.

\bibitem{La} \sc Lamperti, J.W. (1972). \rm Semi-stable Markov processes. {\it Z.
Wahrsch. verw. Gebiete}, {\bf 22}, 205--225.

\bibitem{mi}\textsc{ Millar, P.W.} (1973). Radial processes.\textit{ Annals of Probab.}, \textbf{1},
613-626.


\bibitem{po} \sc Port, S.C. (1969). \rm The First Hitting Distribution of a
Sphere for Symmetric Stable Porcesses. {\it Trans.
Amer. Math. Soc.}, {\bf 135}, 115-125.

\bibitem{ry} \sc  Revuz, D. and Yor, M. (1999). \it Continuous martingales and Brownian motion. \rm Springer-Verlag, Berlin.


\bibitem{Sa} \sc Sato, K.I. (1999). {\it L\'evy processes and infinitely divisible
distributions}. \rm Cambridge University Press, Cambridge.

\bibitem{vi} \sc Vigon, V. (2002).\rm Votre L\'evy rampe-t-il?, {\it J. London Math. Soc.,} {\bf 65}, 243-256.

\end{thebibliography}
\end{document}